\documentclass[titlepage,11pt]{article}
\usepackage{latexsym}
\usepackage{amsfonts}
\usepackage{amsmath}
\newcommand\blackslug{\hbox{\hskip 1pt \vrule width 4pt height 8pt depth 1.5pt
        \hskip 1pt}}
\newcommand\bbox{\hfill \quad \blackslug \bigbreak}
\def\d{\hbox{-}}
\def\c{\hbox{-}\cdots\hbox{-}}
\def\l{,\ldots,}

\title{Colouring perfect graphs with bounded clique number}
\author{Maria Chudnovsky\thanks{Supported by NSF grant DMS-1550991 and US Army Research Office Grant W911NF-16-1-0404.}\\
Princeton University, Princeton, NJ 08544
\\
\\
Aur\'elie Lagoutte\thanks{Partially supported by ANR Grant 
STINT: ANR-13-BS02-0007, and performed while visiting Princeton
University.} \\
LIP, UMR 5668, ENS Lyon, CNRS, UCBL, INRIA, Universit\'e de Lyon, France
\\
\\
Paul Seymour\thanks{Supported by ONR grant N00014-14-1-0084 and 
NSF grant DMS-1265563.}\\
Princeton University, Princeton, NJ 08544
\\
\\
Sophie Spirkl\\
Princeton University, Princeton, NJ 08544}

\date{November 4, 2015; revised September 2, 2016}

\newtheorem{thm}{}[section]

\newcommand{\Proof}{\noindent{\bf Proof.}\ \ }

\begin{document}
\maketitle
\begin{abstract}
A graph is {\em perfect} if the chromatic number of every induced subgraph equals the size of its largest clique, and
an algorithm of Gr\"otschel, Lov\'asz, and Schrijver~\cite{GLS} from 1988 finds an optimal colouring of a perfect graph in polynomial time. But this
algorithm uses the ellipsoid method, and it is a well-known open question to construct a ``combinatorial'' polynomial-time
algorithm that yields an optimal colouring of a perfect graph. 

A {\em skew partition} in $G$ is a partition $(A,B)$ of $V(G)$ such that $G[A]$ is not connected and
$\overline{G}[B]$ is not connected, where $\overline{G}$ denotes the complement graph ; and it is {\em balanced} 
if an additional parity condition of paths in $G$ and $\overline{G}$ is satisfied.

In this paper we first give a polynomial-time algorithm that, with input a perfect graph, 
outputs a balanced skew partition if there is one. Then we use this to obtain a combinatorial algorithm that finds an 
optimal colouring of a perfect graph with clique number $k$, in time that is polynomial for fixed $k$.
\end{abstract}

\section{Introduction}
All graphs in this paper are finite and have no loops or parallel edges. We denote the chromatic number of $G$ by $\chi(G)$, and
the cardinality of the largest clique of $G$ by $\omega(G)$, and a colouring of $G$ with $\omega(G)$ colours is called an {\em optimal} colouring.
If $X\subseteq V(G)$, we denote by $G[X]$
the subgraph of $G$ {\em induced} on $X$, that is, the subgraph with vertex set $X$ and edge set all edges of $G$ with both ends in $X$.
An {\em antipath} in $G$ is an induced path in the complement graph $\overline{G}$.
A {\em hole} in $G$ is an induced cycle of length at least four, and an {\em antihole} in $G$ is a hole in $\overline{G}$.
A graph is {\em Berge} if it has no odd hole or odd antihole (a hole or antihole is {\em odd} if 
it has an odd number of vertices). This definition arose from a conjecture of Berge~\cite{berge}, 
now the strong perfect graph theorem~\cite{CRST}, that every Berge graph is {\em perfect} --- that is, that 
every induced subgraph admits an optimal colouring.

Gr\"otschel, Lov\'asz, and Schrijver~\cite{GLS} showed in 1988 that the ellipsoid method can be applied to find an 
$\omega(G)$-colouring of a perfect graph, and hence of a Berge graph in view of the result of~\cite{CRST}.
Remarkably, however, we have not yet been able to find a ``combinatorial'' algorithm to do the same; and the result of this paper is
a step in that direction. We give:

\begin{thm}\label{mainalg}
An algorithm that, with input an $n$-vertex Berge graph, outputs in time $O(n^{(\omega(G)+1)^2})$ 
an optimal colouring of $G$.
\end{thm}

A {\em skew partition} in $G$ is a partition $(A,B)$ of $V(G)$ such that $G[A]$ is not connected and $\overline{G}[B]$ is not connected.
It is {\em balanced} if in addition:
\begin{itemize}
\item for all nonadjacent $u,v\in B$, every induced path of $G$ with ends $u,v$ and with interior in $A$ has even length, and
\item for all adjacent $u,v\in A$, every antipath of $G$ with ends $u,v$ and with interior in $B$ has even length.
\end{itemize}
(The {\em length} of a path is the number of edges in it.) A skew partition is {\em unbalanced} if it is not balanced.

To construct \ref{mainalg}, we combine an algorithm of Chudnovsky, Trotignon, Trunck and Vu\v{s}kovi\'{c}~\cite{noskew},
that optimally colours a Berge graph with no balanced skew partition, with an algorithm that finds a 
balanced skew partition in a Berge graph if one exists, and an algorithm that combines optimal colourings of two parts of the 
resultant decomposition into an optimal colouring of the whole.  The next six sections are 
devoted to obtaining the latter, and we return to colouring in section 8.

\section{Skew partitions}

Berge's conjecture was proved~\cite{Mjoin,CRST} via a decomposition theorem,
that every Berge graph either admits a balanced skew partition,
or admits one of two other decompositions, or it belongs to one of five well-understood classes. The details of this theorem are not
important here, but one might hope to apply the decomposition theorem to solve other open questions.
Applying the decomposition theorem to algorithmic problems about Berge graphs
has been stalled because until now we have not been able to find a balanced skew partition in 
a Berge graph in polynomial time. Giving such an algorithm is the first main result of this paper.
Some background:
\begin{itemize}
\item there is an algorithm by Kennedy and Reed~\cite{kenreed} that outputs in time $O(n^6)$
a skew partition of a general $n$-vertex graph if one exists (but it is not necessarily balanced);
\item it is NP-hard~\cite{trot} to test if a general graph admits
a balanced skew partition; and
\item there is an algorithm by Trotignon~\cite{trot} and improved by Charbit, Habib, Trotignon, and Vu\v{s}kovi\'{c}~\cite{habib},  
that tests in time $O(n^5)$ whether an $n$-vertex Berge graph admits
a balanced skew partition (but it does not actually output such a partition if one exists).
\end{itemize}

A {\em component} of $G$ is a maximal nonempty subset $X\subseteq V(G)$ such that $G[X]$ is connected, and an {\em anticomponent} of $G$
is a component of $\overline{G}$.
A skew partition $(A,B)$ is {\em loose} if either some vertex in $A$ is adjacent to every vertex 
in some anticomponent of $G[B]$, or some vertex
in $B$ has no neighbours in some component of $G[A]$; and {\em tight} if it is not loose. 
Loose skew partitions can easily be converted to balanced ones in a Berge graph.
(This is explained at the end of section 4.)
We give three algorithms finding skew partitions in this paper, as follows.
\begin{thm}\label{alg1}
An algorithm with input an $n$-vertex Berge graph $G$; in time $O(n^6)$ it outputs a list of all unbalanced tight skew partitions in $G$.
\end{thm}
\begin{thm}\label{alg2}
An algorithm with input an $n$-vertex graph $G$; in time $O(n^6)$ it outputs a list of all tight skew partitions in $G$. (There are at most $n^4$ of them.)
\end{thm}
\begin{thm}\label{alg3}
An algorithm with input an $n$-vertex graph $G$; in time $O(n^6)$ it decides whether $G$ contains a loose skew partition, and if so 
outputs a loose skew partition, balanced if $G$ is Berge.
\end{thm}
In particular, by running \ref{alg3}, and then if necessary comparing the outputs of \ref{alg1} and \ref{alg2}, 
we obtain:

\begin{thm}\label{balancedalg}
An algorithm with input an $n$-vertex Berge graph $G$; in time $O(n^6)$ it outputs a balanced skew partition in $G$, if there is one.
\end{thm}
(Another way to find a balanced skew partition in a Berge graph is, 
after running \ref{alg3}, to just take the skew partitions of the list output by \ref{alg2}, and check directly if any of them are balanced.
This can be done easily, but it seems to take time $O(n^7)$.)
Our methods for \ref{alg2} and \ref{alg3} are both by modifications of the beautiful algorithm of Kennedy and Reed~\cite{kenreed} mentioned above. 

We remark that this shows that a graph has only polynomially-many tight skew partitions. The same is not true for loose skew partitions,
even in a Berge graph, and even if we only count loose balanced skew partitions. For example the star $K_{1,n}$ has exponentially
many loose balanced skew partitions.

In section 7 we sketch how these algorithms can be extended to ``trigraphs'' (graphs with the adjacency of certain pairs of
vertices undecided). 

\section{Unbalanced and tight}

In this section we give the algorithm \ref{alg1}; it outputs all unbalanced, tight skew partitions in a Berge graph.
We begin with some definitions: 
$G$ is {\em anticonnected} if $\overline{G}$ is connected;
if $A,B\subseteq V(G)$, $A$ is
{\em complete} to $B$ if $A\cap B=\emptyset$ and every vertex of $A$ is adjacent to every vertex of $B$,
and $A$ is {\em anticomplete} to $B$ if $A$ is complete to $B$ in $\overline{G}$; 
a vertex $v$ is complete to a set $A$ if $\{v\}$ is complete to $A$, and similarly for anticomplete; and
$N(v)$ or $N_G(v)$ denotes the set of neighbours of a vertex $v$.
We will need the Roussel-Rubio lemma~\cite{R&R}, the following.
\begin{thm}\label{wonderful}
Let $G$ be Berge, let $X\subseteq V(G)$ be anticonnected, and
$P$ be an induced path in $G \setminus X$ with odd length, such that both ends of $P$ are
complete to $X$. Then either:
\begin{itemize}
\item two adjacent vertices of $P$ are complete to $X$; or
\item $P$ has length $\ge 5$ and there are two nonadjacent vertices $a,b\in X$, where $a,b$ is adjacent to the second vertex of $P$ and 
$b$ to the second last, and there are no other edges between $a,b$ and the interior of $P$; or
\item $P$ has length $3$ and there is an odd antipath with interior in $X$, joining the internal vertices of
$P$.
\end{itemize}
\end{thm}

Let us say a skew partition $(A,B)$ in $G$ is {\em square-based} if there is a hole of length four in $G[B]$,
with vertices $a\d b\d c\d d\d a$ in order, such that
\begin{itemize}
\item $B =(N(a)\cap N(c))\cup (N(b)\cap N(d))$; 
that is, every vertex in $B$ is either complete to 
$\{a,c\}$ or complete to $\{b,d\}$, and no vertex in $A$ is complete to either of these sets; and
\item 
there is an odd induced path in $G$ joining $a,c$ with interior in $A$.
\end{itemize}
We need the following theorem.

\begin{thm}\label{squarebased}
Let $(A,B)$ be an unbalanced, tight skew partition in a Berge graph $G$. Then either $(A,B)$ is square-based in $G$, or
$(B,A)$ is square-based in $\overline{G}$.
\end{thm}
\Proof
This is essentially theorem 4.4 of~\cite{CRST}, but we include a proof for the reader's convenience.
Since $(A,B)$ is not balanced, we may assume (passing
to $\overline{G}$ if necessary) that there exist $u,v\in B$ joined by an induced path $P$ with interior in $A$, where $P$
has length at least three and odd. Let $A_1$ be the component of $G[A]$ that contains the interior of $P$; 
and let $B_1$ be the anticomponent of $G[B]$
that contains $u,v$. Let $P$ have vertices $p_1\c p_k$ in order. Let $A_2\ne A_1$ be another component of $G[A]$, and let $B_2\ne B_1$
be another anticomponent of $G[B]$.
\\
\\
(1) {\em If there exist $x,y\in B_2$ such that $x\d p_2\c p_{k-1}\d y$ is an induced path, then $(A,B)$ is square-based.}
\\
\\
For certainly every vertex in $B\setminus B_1$ is complete to $\{u,v\}$, and every vertex in $B_1$ is complete to $\{x,y\}$. 
We must check that no vertex in $A$ is complete to $\{u,v\}$ or to $\{x,y\}$. From the symmetry between $B_1,B_2$ (exchanging $\{u,v\}$
with $\{x,y\}$) it is enough to show that no vertex in $A$ is complete to $\{u,v\}$. Suppose then that $t$ is such a vertex.
Since $t\d p_1\d p_2\c p_k\d t$ is not an odd hole, $t$ has
a neighbour in the interior of $P$, and in particular $t\in A_1$. Since $u,v$
both have neighbours in $A_2$ (because $(A,B)$ is tight), there is an induced path $Q$ joining $u,v$ with interior in $A_2$. Since
$P\cup Q$ is not an odd hole, $Q$ has odd length. But then adding $t$ to $Q$ gives an odd hole, a contradiction. Thus there
is no such $t$, and so $(A,B)$ is square-based. This proves (1).
\\
\\
(2) {\em If $P$ has length at least five then $(A,B)$ is square-based.}
\\
\\
For in this case, since the ends of $P$ are complete to $B_2$, and no internal vertex of $P$ is complete to $B_2$ (because $(A,B)$ is tight),
the Roussel-Rubio lemma~\ref{wonderful} implies that there are nonadjacent $x,y\in B_2$ such that $x\d p_2\c p_{k-1}\d y$
is an induced path, and the claim follows from (1). This proves (2).

\bigskip
Thus we may assume that $P$ has length three, for every choice of $P$; and similarly (passing to $\overline{G}$) that there is no antipath 
of length at least five and odd, with ends in $A$ and interior in $B$. There is an antipath $Q$ joining $p_2,p_3$ with interior in $B_2$,
since $p_2,p_3$ both have a non-neighbour in $B_2$ (because $(A,B)$ is tight); and since $p_2\d p_4\d p_1\d p_3$ is also an antipath,
and its union with $Q$ does not give an odd antihole, it follows that $Q$ has odd length. Consequently $Q$ has length three; let
its vertices be $p_2\d y\d x\d p_3$ in order. But then $x\d p_2\d p_3\d y$ is an induced path, and (1) implies that $(A,B)$
is square-based. This proves \ref{squarebased}.~\bbox

We can generate all tight square-based skew partitions $(A,B)$ in time $O(n^6)$ (where $n=|V(G)|$), as follows.
\begin{itemize}
\item We find all holes $a\d b\d c\d d\d a$ of length four in $G$ (there are at most $n^4$ of them).
\item For each hole $a\d b\d c\d d\d a$, set $B= (N(a)\cap N(c))\cup (N(b)\cap N(d))$ and $A = V(G)\setminus B$, and test whether
$G[A]$ and $\overline{G}[B]$ are both disconnected. (This takes time $O(n^2)$ for each hole.) If not then move on to the next hole.
\item If $G[A]$ and $\overline{G}[B]$ are both disconnected, then $(A,B)$ is a skew partition. 
Now we check whether it is tight. This takes time $O(n^2)$. If not move on to the next hole.
\item If it is tight, we select an induced path $P$ between $a,c$ with interior in $A$ and test whether it is odd, and if so then $(A,B)$
is square-based and we add $(A,B)$
to the output list. This takes time $O(n^2)$. In either case we move on to the next hole.
This is correct, because all induced paths in $G$ between $a,c$ with interior in $A$ have the same parity. To see this,
fix some induced path $P_1$ between $a,c$ with interior in $A_1$, and another, $P_2$,  with interior in $A_2$, where $A_1, A_2$
are distinct components of $G[A]$. (These paths exist since both $a,c$ have neighbours in both $A_1, A_2$, because $(A,B)$ is tight.)
Since $P_1\cup P_2$ is a hole it follows that $P_1,P_2$ have the same parity; and for every other induced path $Q$ joining $a,c$ with 
interior in $A$, one of $Q\cup P_1, Q\cup P_2$ is a hole and so $Q$ has the same parity as both $P_1,P_2$.
\end{itemize}
This procedure generates a list $L$ of all tight, square-based skew partitions in $G$. 
Now we run the same procedure in $\overline{G}$, and for 
each $(C,D)$ in its output list we append $(D,C)$ to $L$, making a list $L'$ say of all skew partitions $(A,B)$ of $G$ that are tight
and either square-based in $G$ or such that $(B,A)$ is square-based in $\overline{G}$.  By \ref{squarebased},
$L'$ is the desired list of all unbalanced tight skew partitions in $G$. This completes the description of \ref{alg1}.

\section{Finding all tight skew partitions}

The algorithms \ref{alg2} and \ref{alg3} both apply the Kennedy-Reed algorithm from~\cite{kenreed}, and next we 
explain these applications. If $X\subseteq V(G)$, we denote $G[V(G)\setminus X]$ by $G\setminus X$.
We say a {\em cutset} of a graph $G$ is a subset $X\subseteq V(G)$ such that $G\setminus X$ is disconnected.
If $(A,B)$ is a skew partition of $G$, and 
there are vertices $a_1,a_2\in A$, in different components of $G[A]$, 
and an anticomponent $B_1$ of $G[B]$, such that both $a_1,a_2$ are complete to $B_1$, then $B$ is called a {\em T-cutset} in $G$.

The Kennedy-Reed algorithm~\cite{kenreed} generates a list $L$ of subsets of the vertex set of an $n$-vertex graph $G$, with the following properties: 
\begin{itemize}
\item $L$ has at most $n^4$ members;
\item each member $X$ of $L$ is a cutset of $G$;
\item for every skew partition $(A,B)$ of $G$ such that $B$ is not a T-cutset, some member of $L$ is a subset of $B$.
\end{itemize}
It can be applied to any graph, not only to Berge graphs; and its running time is $O(n^6)$. (Incidentally,
Kennedy and Reed erroneously claimed 
a running time of $O(mn^4)$ where $m$ is the number of edges.)
We describe it later; but let us see first how it can be used for \ref{alg2} and \ref{alg3}.

For \ref{alg2} we must generate a list of all tight skew partitions. We use:
\begin{thm}\label{tight}
Let $(A,B)$ be a tight skew partition in a graph $G$. Then no proper subset of $B$ is a cutset of $G$.
\end{thm}
\Proof Suppose that $X\subseteq B$, and there exists $b\in B\setminus X$. We claim that $G\setminus X$ is connected.
For let $A_1$ be a component of $G[A]$, and let $C$ be the component of $G\setminus X$ that includes $A_1$. 
Since every vertex in $B\setminus X$ has a neighbour in $A_1$ (because $(A,B)$ is tight) it follows that 
$B\setminus X\subseteq C$. For every component $A_2$ of $G[A]$, since $b\in C$ has a neighbour in $A_2$ (because
$(A,B)$ is tight), it follows that $A_2\subseteq C$. Consequently $G\setminus X = C$, and so $X$ is not a cutset. 
This proves \ref{tight}.~\bbox

We can generate a list of all tight skew partitions in $G$ as follows. First, we observe that if $(A,B)$ is a tight skew partition in $G$
then $B$ belongs to the list $L$ output by the Kennedy-Reed algorithm; because $B$ is not a T-cutset, so
some subset of $B$ belongs to $L$,
by the properties of $L$, and no proper subset is a cutset by \ref{tight}. Thus, all we have to do is check which
members $B$ of the list $L$ give rise to tight skew partitions; that is, we must check, for each $B\in L$, whether
$G[B]$ is not anticonnected and $(A,B)$ is tight, setting $A=V(G)\setminus B$.
This takes time $O(n^2)$ for each choice of $B$, and there are at most $n^4$ choices of $B$ to check.

This completes the description of \ref{alg2}. We note that it can be applied to any graph, not only to Berge graphs.

It follows that there are at most $n^4$ tight skew partitions in any $n$-vertex graph, because the corresponding cutset of each belongs to
the list output by the Kennedy-Reed algorithm. We can improve this. While the list of the Kennedy-Reed algorithm 
can indeed have $O(n^4)$ terms, there are in fact
at most $n^3\log n$ tight skew partitions; because for every skew partition $(A,B)$, $B$ appears at least $k_2$ times
in the list output by Kennedy-Reed, where $k_2$ is the size of the second largest anticomponent of $G[B]$. (We omit further details.)
We do not know how many tight skew partitions there can really be, and in fact the authors do not know of a graph with more than linearly many 
tight skew partitions.

\section{Finding a loose skew partition}

A {\em star cutset} in $G$ is a cutset $B$ with $|B|\ge 2$ such that some vertex $b$ in $B$ is adjacent to all other vertices in $B$.
If $B$ is a star cutset, then $(A,B)$ is a skew partition where $A=V(G)\setminus B$, and $(A,B)$ is loose since
$\{b\}$ is an anticomponent of $B$, and either some vertex in $A$ is complete to $\{b\}$, or $b$ is anticomplete to $A$.
If $v$ is a vertex, $N(v)$ denotes its set of neighbours, and 
$N[v]=N(v)\cup \{v\}$. 
For \ref{alg3}, it is helpful first to find a star cutset if there is one; we begin with that. The following is due to 
Chv\'atal~\cite{chvatal}.

\begin{thm}\label{starcutalg}
An algorithm that, with input an $n$-vertex graph $G$, outputs in time $O(n^3)$ a star cutset of $G$ if there is one.
Moreover, the same algorithm, with input a Berge graph $G$, outputs a star cutset $B$ of $G$ such that $(A,B)$ is balanced,
if there is a star cutset.
\end{thm}
Here is the algorithm. For each vertex $v$ of $G$ in turn, we do the following:
\begin{itemize}
\item if $N(v)\ne \emptyset$ and $G\setminus N[v]$ is not connected, output $N[v]$ and stop;
\item if $|N(v)|\ge 2$, and $G\setminus N[v]$ has a unique component $C$ say, and some vertex $u\in N(v)$ is anticomplete to $C$,
output $N[v]\setminus \{u\}$ and stop;
\item if $|N(v)|\ge 3$, and $N[v] = V(G)$, and there are two nonadjacent vertices $x,y\in N(v)$, output $V(G)\setminus \{x,y\}$ and stop;
\item Otherwise move on to the next vertex.
\end{itemize}
If we have examined all vertices with no output, return that there is no star cutset.

This completes the description of \ref{starcutalg}. It is easy to see that it works in a general graph. To see the second claim,
that in a Berge graph
it returns a balanced star cutset, let $B$ be the star cutset returned, let $v$ be the vertex that was being examined at that stage,
and let $A=V(G)\setminus B$.
Let $A_1\l A_r$ be the components of $G[A]$, and let $B_1\l B_s$ be the anticomponents of $G[B]$, where $B_1 = \{v\}$.
From the description of the algorithm, every neighbour of $v$ in $A$ is the vertex of a singleton component of $G[A]$, and therefore $v$ 
is anticomplete to every component of $G[A]$ with more than one vertex.
To check that $(A,B)$ is balanced, we have to check that there is no odd induced path $P$ with ends in $B$ and interior in $A$ with length
more than one; and there is no odd antipath $Q$ with ends in $A$ and interior in $B$ with length more than one. Suppose that
some such $P$ or $Q$ exists. If $P$ exists, then both ends of $P$ belong to the same anticomponent of $B$, and its interior belongs 
to one component of $G[A]$; and similarly if $Q$ exists then its ends belong to the same component of $G[A]$, and its interior to
one anticomponent of $G[B]$. Thus we may choose $A_i$ and $B_j$ such that $A_i\cup B_j$ includes $P$ or $Q$. Now $|B_j|>1$, so 
$j\ne 1$; and $|A_i|>1$, so $v$ has no neighbour in $A_i$. But then adding $v$ to $P$ or $Q$ gives an odd hole or antihole, a contradiction.
This proves the second claim of \ref{starcutalg}.

\bigskip

We also need a result of~\cite{hoang}:
\begin{thm}\label{Tcutalg}
An algorithm that with input an $n$-vertex graph, outputs in time $O(n^5)$ a T-cutset of $G$, if there is one.
\end{thm}
Here is the algorithm.
List all the nonadjacent pairs of distinct vertices of $G$. For each such pair $(a_1,a_2)$,
list the anticomponents of $G[N(a_1)\cap N(a_2)]$. For each such anticomponent $B_1$ say,
let $B_2$ be the set of all vertices in $V(G)\setminus (B_1\cup \{a_1,a_2\})$ that are
complete to $B_1$. If either $B_2=\emptyset$, or there is a path in $G\setminus (B_1\cup B_2)$ between $a_1, a_2$,
move on to the next anticomponent. Otherwise
output $B_1\cup B_2$ and stop. It is easy to see that this performs as claimed. 

\bigskip
We remark that it is not true that this algorithm always returns a balanced skew partition in a Berge graph with a T-cutset.
Let us return to finding a loose skew partition. We use the following:
\begin{thm}\label{findloose}
Let $G$ be a connected $n$-vertex graph with no star cutset, and let $L$ be the output of the Kennedy-Reed algorithm applied to $G$. 
If there is a loose skew partition
in $G$, then 
$(A,B)$ is a loose skew partition for some $B\in L$, where $A=V(G)\setminus B$.
\end{thm}
\Proof
Let $(A,B)$ be a loose skew partition in $G$. From the properties of $L$, there is a cutset $X\in L$ with $X\subseteq B$.
Suppose for a contradiction that $G[X]$ is anticonnected.  Let $B_1$ be an anticomponent of $B$ with $X\subseteq B_1$.
Now $X\ne \emptyset$ since $X$ is a cutset and $G$ is connected. Let $v\in B\setminus B_1$, chosen with a neighbour in $A$ if possible.
Since there is no star cutset in $G$ it follows that $X\cup \{v\}$ is not a star cutset, and so 
$G\setminus (X\cup \{v\})$ is connected. 
Now $G\setminus X$ is not connected since $X$ is a cutset; and so $G\setminus X$
has exactly two components, and one of them consists of the singleton $v$. In particular, $v$ is anticomplete to $A$ (and so from
the choice of $v$, $B\setminus B_1$ is anticomplete to $A$); and also $v$ is anticomplete to $B_1\setminus X$. Since
$v$ is complete to $B_1$ we deduce that $X=B_1$. But then every component of $G\setminus B$ is a component of $G\setminus X$, and
there are at least two such components, contradicting that $G\setminus X$ has only one component not containing $v$.

This proves that $G[X]$ is not anticonnected. Consequently $(V(G)\setminus X, X)$ is a skew partition, and we claim it is 
loose. For let $X_1$ be
an anticomponent of $X$, and let $B_1$ be the anticomponent of $B$ with $X_1\subseteq B_1$. If there is a vertex
$v\in B\setminus (B_1\cup X)$, then $v\notin X$, and $v$ is complete to an anticomponent of $X$, showing that $(V(G)\setminus X, X)$
is loose. Thus we may assume that $B\setminus (B_1\cup X)=\emptyset$, that is, $B\setminus B_1\subseteq X$.
In particular, there is a second anticomponent $X_2$ of $X$
with $X_2\subseteq B\setminus B_1$; and by the same argument with $X_1$ replaced by $X_2$, we may assume that $B_1\subseteq X$, that is,
$B=X$. But then $(V(G)\setminus X, X)$ is loose since $(A,B)$ is loose. This proves our claim that $(V(G)\setminus X, X)$ is loose;
and hence the theorem holds.  This proves \ref{findloose}.~\bbox

Thus we have:
\begin{thm}\label{loosealg}
An algorithm that, with input a connected $n$-vertex graph $G$ with no star cutset, 
outputs in time $O(n^6)$ a loose skew partition of $G$ if one exists.
\end{thm}
The algorithm is as follows. Use \ref{Tcutalg} to test if $G$ has a T-cutset, and if \ref{Tcutalg} returns a T-cutset $B$,
output $(V(G)\setminus B)$ and stop. Otherwise,
run the Kennedy-Reed algorithm, and let $L$ be the list it outputs. For each $B\in L$, let $A=V(G)\setminus B$,
and test if $(A,B)$ is loose. If so, output it and stop. If we exhaust $L$ with no output, return that there is no loose skew partition.

\bigskip

Every star cutset gives a loose skew partition; because if $B$ is a star cutset in $G$, and $\{v\}$ is an anticomponent of $G[B]$,
either some vertex in $A$ is complete to $\{v\}$, or $v$ is anticomplete to $A$. Moreover, every disconnected graph with at least
one edge and at least four vertices has a star cutset. Thus combining \ref{starcutalg} and \ref{loosealg}
gives:
\begin{thm}\label{loosealg2} An algorithm that, with input an $n$-vertex graph $G$,
outputs in time $O(n^6)$ a loose skew partition of $G$ if one exists.
\end{thm}

To complete \ref{alg3}, we must output a loose balanced skew partition when the input graph is Berge and has a loose skew partition. 
We do so as follows.
First we run \ref{starcutalg}, both in $G$ and in $\overline{G}$, and if either returns a star cutset we return the corresponding
balanced skew partition and stop. Thus we may assume that neither of $G,\overline{G}$ has a star cutset. We may also assume that
$G$ has at least one edge and at least four vertices; and since $G$ has no star cutset, it follows that 
$G$ is connected and similarly so is $\overline{G}$. Now we run \ref{loosealg}, and we may assume it returns a loose skew partition 
$(A,B)$. By passing to the complement if necessary, we may assume that for some anticomponent $B_1$ of $G[B]$, and some component
$A_1$ of $G[A]$, some vertex in $B_1$ is anticomplete to $A_1$. Now do the following:
\begin{itemize}
\item If some vertex $v\in B\setminus B_1$ has no neighbours in $A_1$, replace $A$ by $A'=A\cup \{v\}$ and $B$ by $B'=B\setminus \{v\}$.
Then $(A',B')$ is a skew partition (note that $B\ne B_1\cup \{v\}$ since $G$ has no star cutset); and $B_1$ is an anticomponent 
of $B'$, and $A_1$ is a component of $A'$.
\item If some vertex $v\in A$ is complete to some anticomponent $B_2\ne B_1$ of $B$, and not complete to $B_1$, replace
$A$ by $A' = A\setminus \{v\}$ and $B$ by $B\cup \{v\}$. 
Then $(A',B')$ is a skew partition (note that $G[A']$ has at least two components since $\overline{G}$ has no star cutset); and $B_1$
is contained in an anticomponent $B_1'$ of $G[B']$; and some vertex of $B_1'$ has no neighbour in any of the (at least one) components
of $G[A']$ that are included in $A_1$.
\end{itemize}
We repeat this process until there is no further progress; it repeats at most $2n$ times, since at each step the quantity
$2|B_1|-|B|$ is increased by 1. Consequently this takes time $O(n^3)$. Let it terminate with a skew partition $(A,B)$; then
by theorem 4.2 of~\cite{CRST}, $(A,B)$ is balanced. (Not quite; theorem 4.2 of~\cite{CRST} assumes that we have chosen $(A,B)$ with $2|B_1|-|B|$ maximum, but all the proof in that paper needs is 
that neither of the two operations above can be applied to increase $2|B_1|-|B|$.) Moreover $(A,B)$ is loose, from the way we found it; we output it and stop.
This completes \ref{alg3}.

\section{The Kennedy-Reed algorithm}
In this section we describe the Kennedy-Reed algorithm, which we have slightly modified.
A {\em clique cutset} of $G$ is a clique $X$ of $G$ that is a cutset.
We need an algorithm of Tarjan, and first we need some definitions to describe that.
A {\em CC-decomposition tree} of a graph $G$ (CC stands for clique cutset) consists of
a rooted tree $T$, and for each $v\in V(T)$, a subset $X_v$ of $V(G)$, satisfying the following conditions:
\begin{itemize}
\item every vertex of $T$ has either two or no children (the {\em children} of a vertex $v$ 
are the neighbours of $v$ that do not lie on the path between $v$ and the root);
\item if $w$ is the root of $T$ then $X_w = V(G)$, and if $u$ is a child of $v\in V(T)$ then $X_u$ is a proper subset of $X_v$;
\item for each $s\in V(T)$ with children $r,t$, $X_r\cup X_t = X_s$, and $X_r\cap X_t$ is a clique cutset of $G[X_s]$, and 
$X_r\setminus X_t$ is anticomplete to $X_t\setminus X_r$;
\item if $v\in V(T)$ has no children, then $G[X_v]$ has no clique cutset.
\end{itemize}

Tarjan~\cite{tarjan} gave the following:
\begin{thm}\label{tarjanalg}
An algorithm that, given a graph $G$, outputs in time $O(|V(G)|\cdot |E(G)|)$ a CC-decomposition tree such that at most
$|V(G)|-2$
vertices of the tree have children.
\end{thm}

We need the following two lemmas:
\begin{thm}\label{cutsetgrow}
Let $T$ and $(X_t:t\in V(T))$ form a CC-decomposition tree for $G$. For each $t\in V(T)$, every clique cutset of $G[X_t]$ is a clique cutset of $G$.
\end{thm}
\Proof
It suffices to prove that if $t$ is a child of $s$ and $Z$ is a clique cutset of $G[X_t]$ then $Z$ is a clique cutset of $G[X_s]$.
Let the children of $s$ be $r,t$, and let $Y = X_r\cap X_t$. Since $Y$ is a clique, there is a component $A_1$ of $G[X_t]\setminus Z$ that 
includes $Y\setminus Z$; let $A_2 = X_t\setminus (A_1\cup Z)$. Thus $A_2\ne \emptyset$ (since $Z$ is a clique cutset of $G[X_t]$), and since 
$A_1\cup (X_r\setminus X_t)$ is anticomplete to $A_2$, it follows that $Z$ is a clique cutset of $G[X_s]$. This proves \ref{cutsetgrow}.~\bbox

\begin{thm}\label{decomptree}
Let $T$ and $(X_t:t\in V(T))$ form a CC-decomposition tree for $G$, and let $Z$ be a clique cutset of $G$. Then there exists
$s\in V(T)$, with children $r,t$ say, such that $X_r\cap X_t\subseteq Z$.
\end{thm}
\Proof 
We may assume that $Z$ is a minimal clique cutset of $G$.
Choose $s\in V(T)$ with $X_s$ minimal such that $Z$ is a clique cutset of $G[X_s]$.
Since $Z$ is a minimal clique cutset of $G$, \ref{cutsetgrow} implies that $Z$ is a minimal clique cutset of $G[X_s]$. 
Since $G[X_s]$ has a clique cutset, it follows that $s$ has children $r,t$ say. 
Let $A_1\l A_k$ be the components of $G[X_s\setminus Z]$; thus $k\ge 2$.
We assume for a contradiction that $X_r\cap X_t\not \subseteq Z$. Since $X_r\cap X_t\setminus Z$
is a clique, all its vertices belong to the same one of $A_1\l A_k$, say $A_1$. 
Since $Z$ is a clique and 
$Z\subseteq X_r\cup X_t$, and $X_r\setminus X_t$ is anticomplete to $X_t\setminus X_r$, it follows that $Z$ is a subset of 
one of $X_r, X_t$, say $X_r$. Since $X_r\cap X_t\setminus Z$
is nonempty, it follows that 
$X_r\cap A_1\ne \emptyset$. Since $Z$ is not
a clique cutset of $G[X_r]$ from the choice of $s$, it follows that $A_2\l A_k$ are disjoint from $X_r$. 
Since $k\ge 2$, the argument with $X_r, X_t$ exchanged implies that $Z\not\subseteq X_t$, and so there exists $v\in Z\setminus X_t$.
But $v$ has a neighbour in $A_2$, since $k\ge 2$ and $Z$ is a minimal clique cutset of $G[X_s]$. This is impossible
since $X_r\setminus X_t$ is anticomplete to $X_t\setminus X_r$. 
Consequently $X_r\cap X_t \subseteq Z$. This proves \ref{decomptree}.~\bbox

From \ref{tarjanalg}, \ref{cutsetgrow} and \ref{decomptree} we obtain:
\begin{thm}\label{tarjanalg2}
An algorithm that, with input an $n$-vertex graph $G$, outputs in time $O(n^3)$ a list of at most $n-2$ clique cutsets of $G$, 
such that every clique cutset of $G$ includes one of them.
\end{thm}
To see this, we take the output of \ref{tarjanalg}, say $T$ and $(X_t:t\in V(T))$; for each $s\in V(T)$ with children $r,t$ say,
let $K_s=X_r\cap X_t$; and we output the list of all such sets $K_s$. This is correct by \ref{decomptree} and the fact
that only $n-2$ vertices in $T$ have children.

Now we turn to the idea of Kennedy and Reed. Let $G$ be a graph, let $r\in V(G)$, and let $1\le k_2\le k_1\le n$ be integers.
Let $H(k_1,k_2,r)$ be the graph with vertex set $V(G)$, in which a pair $u,v$ of distinct vertices is adjacent if either
\begin{itemize}
\item $u,v$ are adjacent in $G$, or
\item some anticomponent of $G[N(u)\cap N(v)]$ has cardinality at least $k_1$, or
\item some anticomponent of $G[N(u)\cap N(v)]$ has cardinality at least $k_2$ and contains $r$.
\end{itemize}
Now let $(A,B)$ be a skew partition in $G$, such that $B$ is not a T-cutset. 
If the largest anticomponent of $G[B]$ has cardinality $k_1$, and the second largest has
cardinality $k_2$, and the latter contains $r$, it is easy to see that $B$ is a clique cutset of $H(k_1,k_2,r)$.
Let us run \ref{tarjanalg2} on $H(k_1,k_2,r)$ for all choices of $k_1,k_2,r$, and take the union $L$ of their
output lists. Every member of $L$ is a cutset of  $H(k_1,k_2,r)$ for some choice of $k_1,k_2,r$, and hence is a cutset of $G$
since every edge of $G$ is an edge of $H(k_1,k_2,r)$. Moreover, there are at most 
$n^3$ choices of $k_1,k_2,r$, and each choice gives a list via \ref{tarjanalg2} of cardinality at most $n$, so
$L$ has cardinality at most $n^4$. 
Finally, for every skew partition $(A,B)$ such that $B$ is not a T-cutset, there is a choice of $k_1,k_2,r$ such that
$B$ is a clique cutset of $H(k_1,k_2,r)$, and therefore $B$ includes a member of $L$.
Consequently $L$ has the properties described at the start of section 4.

\section{Trigraphs}

A {\em trigraph} consists of a set $V(G)$ of {\em vertices}, and a classification of each pair of distinct vertices
as strongly adjacent, strongly nonadjacent, or semiadjacent. A {\em realization} of a trigraph $G$ is a graph $H$ with $V(H)=V(G)$,
such that every strongly adjacent pair of vertices in $G$ is adjacent in $H$, and every strongly nonadjacent pair in $G$
is nonadjacent in $H$. (Think of the semiadjacent pairs as pairs whose adjacency is undecided; a realization results from
making the decision, for each undecided pair.)

Trigraphs were introduced by one of us in~\cite{Maria}, because of difficulties that arose from the decomposition theorem 
for perfect graphs. In particular, for some of the decompositions used in~\cite{CRST}, the graph is best viewed as being decomposed
into smaller trigraphs rather than into smaller graphs; and trigraphs naturally arise in the context of decomposition theorems
for Berge graphs. Thus no doubt we will eventually need an algorithm to find a balanced skew partition in a Berge trigraph;
and so let us sketch here what needs to be modified to make the algorithms of this paper work for trigraphs.

First, a trigraph is {\em Berge} if every realization is Berge. 
A trigraph $G$ is {\em connected} if there is no partition of $V(G)$
into two nonempty sets $A,B$ such that every vertex in $A$ is strongly nonadjacent to every vertex in $B$.
A {\em component} $X$ of $G$ is a maximal non-null subset of 
$V(G)$ such that $G[X]$ is connected. A vertex $v$ is {\em strongly complete}
to a set $X$ if $v\notin X$ and $v$ is strongly adjacent to every member of $X$;
and $v$ is {\em weakly complete} to $X$ if $v\notin X$ and no vertex in $X$ is strongly antiadjacent to $v$.
The complement $\overline{G}$ of a trigraph $G$ is a second trigraph defined in the natural way, keeping
the set of semiedges unchanged.
A set is {\em anticonnected} if it connected in $\overline{G}$, and {\em anticomponent, strongly anticomplete, weakly anticomplete}
are defined similarly.
A {\em skew partition} of a trigraph $G$
is a partition $(A,B)$ of $V(G)$ such that $G[A],\overline{G}[B]$ are not connected. (Consequently it is a skew
partition in every realization.) It is {\em balanced} if 
it is balanced in every realization. 
A skew partition $(A,B)$ is {\em loose} if either some vertex in $A$ is weakly complete to some anticomponent of $G[B]$, or some vertex in $B$
is weakly anticomplete to some component of $G[A]$; and {\em tight} otherwise.

\ref{squarebased} still holds, with basically the same proof, requiring that the ``hole of length four'' in the definition of 
square-based be formed by strongly adjacent pairs,
and by interpreting ``induced path'' to mean ``induced path in some realization''.
A step of that proof needs the Roussel-Rubio lemma~\ref{wonderful}, but we do not need to extend the Roussel-Rubio lemma to trigraphs; we apply 
it instead to the realization $H$ in which all semiadjacent pairs of $G$ are nonadjacent except for those in the path in question, 
which we make adjacent. Consequently \ref{alg1} still works. 

A {\em cutset} is a set $B\subseteq V(G)$ such that $G\setminus B$ is not connected;
and a {\em T-cutset} is a cutset $B$ with an anticomponent $B_1$ and two vertices $a_1,a_2$ in different components of $G[A]$,
such that $a_1,a_2$ are strongly complete to $B_1$. 

To run the Kennedy-Reed algorithm on a trigraph $G$, we define the graphs $H(k_1,k_2,r)$ (they are still graphs, not trigraphs) 
by saying that $u,v$ are adjacent if either they are
strongly adjacent or semiadjacent in $G$, or their common strong neighbours have an anticomponent of size $\ge k_1$, or one of size $\ge k_2$ 
containing $r$. Otherwise Kennedy-Reed runs as for graphs.

To run \ref{alg2}, we observe that \ref{tight} still holds, with 
the same proof, and \ref{alg2} runs as before. Let us turn to \ref{alg3}.
A {\em star cutset} is defined as before, and \ref{starcutalg} works as before. 
Also \ref{Tcutalg}, \ref{findloose}, and \ref{loosealg} work as before.
To complete \ref{alg3} we also need an algorithm producing a loose balanced skew partition from a loose skew partition 
in a Berge trigraph, and for this we need a trigraph version of
theorem 4.2 of~\cite{CRST}; it does not seem enough to apply the graph version of that theorem to some appropriate realization as
far as we can see. However, such a trigraph version is true, and was proved in~\cite{Maria} (the proof for graphs works
for trigraphs with very little adjustment). Thus all the algorithms finding skew partitions in this paper 
can be extended to trigraphs with little or no work, and with the same running times.

Incidentally, it is tempting to redefine ``loose'' using strong completeness and anticompleteness rather than weak, and ask if we can
still generate the tight skew partitions (now there are more of them). We can, but it needs more care; for instance, star cutsets
are no longer necessarily loose. We omit further details.

\section{Return to colouring}

Now we show how to use \ref{balancedalg} to obtain \ref{mainalg}. We combine it with an algorithm of 
Chudnovsky, Trotignon, Trunck and Vu\v{s}kovi\'{c}~\cite{noskew}, the following:
\begin{thm}\label{noskewalg}
An algorithm that, with input an $n$-vertex Berge graph $G$ with no balanced skew partition, outputs in time $O(n^7)$ an 
optimal colouring of $G$.
\end{thm}

Given a Berge graph that we need to colour, we first apply \ref{balancedalg}, and if it tells us that there is no balanced skew partition,
we just apply \ref{noskewalg} and we are done. The question is, what do we do if \ref{balancedalg} gives us a skew partition $(A,B)$?
In this case, we partition $A$ into two nonempty sets $A_1,A_2$ such that there are no edges between $A_1$ and $A_2$,
and work with the subgraphs induced on $A_1\cup B$ and $A_2\cup B$. With some sort of induction, we obtain optimal colourings
of these two subgraphs; and then we need to piece them together to obtain a colouring of the original graph $G$. So, there are two 
issues: given colourings of these two subgraphs,
how do we fit them together? And
how can we arrange an induction that will lead to a polynomial-time algorithm? 
(At the moment we do not
know how to handle either issue in time polynomial in both $|V(G)|$ and $\omega(G)$; our algorithms are just 
polynomial-time
in $|V(G)|$ for $\omega(G)$ fixed.) We handle the two issues
in sections 9 and 10 respectively.

\section{Making two colourings match}

Fix $k\ge 0$. We assume that we have an
algorithm $\mathcal{P}$ say, that will optimally colour any $n$-vertex Berge graph with clique number at most $k-1$, 
in time $O(n^c)$, where $c\ge k$. Let $(A,B)$ be a balanced skew partition of an $n$-vertex  Berge graph $G$ with clique number $k$.
Let $(A_1, A_2)$ be a partition of $A$ with $A_1, A_2\ne \emptyset$, such that $A_1$ is anticomplete to $A_2$, and for
$i = 1,2$ let $G_i = G[A_i\cup B]$. (Thus $k\ge 2$.)
For $i = 1,2$ let $\phi_i$ be an optimal colouring of $G_i$, mapping from $V(G_i)$ to the set of integers $\{1\l \omega(G_i)\}$. 
Let us describe how to obtain an optimal colouring of $G$.
\begin{enumerate}
\item Take a partition $(B_1, B_2)$ of $B$ into two nonempty sets, such that $B_1$ is complete to $B_2$,
and compute $b_i=\omega(G[B_i])$ for $i = 1,2$. Let $b_i = \omega(G[B_i])$. Exchange $B_1,B_2$ if necessary 
to arrange that $|B_1|-b_1\le |B_2|-b_2$.
\item For $i = 1,2$, let $L_i = \{\phi_i(v):v\in B_1\}$, and let $\ell_i = |L_i|$; 
by permuting the colours, arrange that $L_i = \{1\l \ell_i\}$. 
Let $S_i$ be the set of $v\in A_i\cup B$ with $\phi_i(v) \in L_i$.
So $B_1\subseteq S_i$ and $B_2\cap S_i=\emptyset$. 
\item For $i = 1,2$, let $H_i$ be obtained from $G_i[S_i]$ by adding $\ell_i-b_1$ new vertices,
each complete to $B_1$ and anticomplete to $S_i\setminus B_1$, and all adjacent to one another. Apply $\mathcal{P}$ to $H_i$ to obtain an $\ell_i$-colouring $\xi_i$ of $H_i$.
By permuting colours arrange that the $\ell_i-b_1$ new vertices have colours $b_1+1\l \ell_i$. For each $v\in A_i\cup B$, 
let $\psi_i(v) = \phi_i(v)$
if $v\notin S_i$, and $\psi_i(v) = \xi_i(v)$ if $v\in S_i$. 
\item For $i = 1,2$, let $T_i$ be the set of vertices $v\in V(G_i)$ with $\psi_i(v)\in \{1\l b_1\}$. Apply $\mathcal{P}$
to $G[T_1\cup T_2]$ and to $G\setminus (T_1\cup T_2)$, to obtain a $b_1$-colouring of $G[T_1\cup T_2]$ and a $(k-b_1)$-colouring of
$G\setminus (T_1\cup T_2)$. Combine them to make a $k$-colouring of $G$.
\end{enumerate}

Let us fill in some more detail and explanation. In step 1, such a partition exists because $G[B]$ is not anticonnected, and 
we can find one in time $O(n^2)$. Computing $\omega(G[B_i])$ for $i = 1,2$ takes time at most $O(n^k)$, by trying all subsets of size at most $k$.

In step 2, $S_i$ is the set of vertices of $G_i$ that have the same colour under $\phi_i$ as some vertex in $B_1$. Since $B_2$
is complete to $B_1$, no vertex in $B_2$ belongs to $S_i$.

In step 3, $H_i$ is Berge, since $(A,B)$ is balanced. Also, $\omega(H_i)\le \ell_i$, since every clique of $H_i$ that contains no new vertex
has already been coloured with $\ell_i$ colours, and every clique that contains a new vertex is disjoint from $A\cup B_2$ and has at most $b_1$ vertices in $B_1$.
Also $\ell_i\le k-1$,  because $B_2\ne\emptyset$ and no colour appears under $\phi_i$ in both $B_1, B_2$. Consequently $\mathcal{P}$
can be applied to $H_i$, and it yields an $\ell_i$-colouring. We claim that $|V(H_i)|\le n$. There are only $\ell_i-b_1$ vertices
of $H_i$ that are not vertices of $G$, and there are at least $|B_2|$ vertices in $G$ that are not in $H_i$; and 
$$\ell_i-b_1\le |B_1|-b_1\le |B_2|-b_2\le |B_2|.$$ 
This proves that $|V(H_i)|\le n$, and so this application of $\mathcal{P}$ takes time $O(n^c)$.
Since the new vertices are a clique and so all have different colours,
we can arrange by permuting colours that the $\ell_i-b_1$ new vertices have colours $b_1+1\l \ell_i$. Consequently, only colours $1\l b_1$
appear in $B_1$. Since only colours $1\l \ell_i$ appear in $S_i$ under $\xi_i$, and only colours $\ell_i+1\l k$ appear in $V(G_i)\setminus S_i$
under $\phi_i$, it follows that $\psi_i$ defined in step 3 is a $k$-colouring of $G_i$, and under it only colours $1\l b_1$ appear in $B_1$.
This step takes time $O(n^c)$.

In step 4, since $G_i[T_i]$ is coloured with only $b_1$ colours under $\psi_i$, its clique number is at most $b_1$; and since every clique included in $T_1\cup T_2$
is included in one of $T_1,T_2$, it follows that 
the clique number of $G[T_1\cup T_2]$ is at most $b_1$ (in fact, exactly $b_1$). Similarly, the clique number of $G\setminus (T_1\cup T_2)$
is at most $k-b_1$. Since $b_1, k-b_1<k$, these applications of $\mathcal{P}$ are valid, and yield colourings as described. This step 
takes time $O(n^c)$. 

Consequently, the whole algorithm takes time $O(n^c)$, since $c\ge k\ge 2$.

\section{The induction}

We need the following lemma, a result of Duchet and Meyniel~\cite{duchet} (we give a proof for the reader's convenience, and since it is not explicitly proved in~\cite{duchet}).
The {\em stability number} of $G$ is $\omega(\overline{G})$.
\begin{thm}\label{growstable}
Let $G$ be a nonnull connected graph with stability number $\alpha$. Then there is a connected induced subgraph with at most $2\alpha-1$ vertices
and with stability number $\alpha$.
\end{thm}
\Proof
Choose a stable set $S$ with the following properties:
\begin{itemize}
\item there is a subset $T\subseteq V(G)$ with $|T|<|S|$ such that $G[S\cup T]$ is connected;
\item there is a stable set of cardinality $\alpha$ including $S$;
\item $S$ is maximal subject to these two conditions.
\end{itemize}
This is possible because setting $S$ to be a singleton subset of some largest stable set satisfies the first two bullets.
We claim that $|S|=\alpha$; for suppose not, and let $S'$ be a stable set of cardinality $\alpha$ including $S$. Since 
$S'$ is a maximum stable set, every vertex of $G$ either belongs to $S'$ or has a neighbour in $S'$. 
Since $G$ is connected, there is a shortest path between $S$ and $S'\setminus S$, with vertices $p_1\c p_k$ in order say,
where $p_1\in S$ and $p_k\in S'\setminus S$. Since $S'$ is stable, $k\ge 3$; and since every vertex of $G$ belongs to $S'$
or has a neighbour in $S'$, $k\le 4$. If $k=3$, we could replace $S$ by $S\cup \{p_3\}$ and $T$ by $T\cup \{p_2\}$, 
contrary to the choice of $S$. So $k=4$. If $p_3$ has more than one neighbour in $S'\setminus S$, say $p_4$ and $p_4'$, then
we could replace $S$ by $S\cup \{p_4,p_4'\}$ and $T$ by $T\cup \{p_2,p_3\}$, again contrary to the maximality of $S$. So
$p_4$ is the only neighbour of $p_3$ in $S'\setminus S$. But then $(S'\setminus \{p_4\})\cup \{p_3\}$
is a stable set of cardinality $\alpha$, and we could replace $S$ by $S\cup \{p_3\}$ and $T$ by $T\cup \{p_2\}$, again a contradiction.
This proves that $|S|=\alpha$, and so proves \ref{growstable}.~\bbox

Let us say a {\em $k$-pellet} in a graph $G$ is a subset $P\subseteq V(G)$ with $|P|=2k$ such that $G[P]$ is anticonnected 
and $\omega(G[P])\ge k$. By \ref{growstable} applied in the complement, every anticonnected graph $G$ 
with at least $2\omega(G)$ vertices has
at least one $\omega(G)$-pellet; and only polynomially many (for fixed $\omega(G)$), since $\omega(G)$-pellets have bounded size.
Essentially, we are going to prove inductively that the running time of our algorithm is at most proportional to the number of 
$\omega(G)$-pellets in $G$ (actually, proportional to that number plus one). But to make the argument clearer, let us
replace the inductive step by the following tree structure.

Let us say a {\em SP-decomposition tree} of $G$ (SP stands for skew partition) is a rooted tree $T$, together with a choice of a subset $X_t\subseteq V(G)$
for each $t\in V(T)$, satisfying the following conditions:
\begin{itemize}
\item every vertex of $T$ has two or zero children (the {\em children} of a vertex $v$ are its neighbours that are not on the path between
$v$ and the root);
\item $X_w=V(G)$ if $w$ is the root, and for all $s\in V(T)$, if $r$ is a child of $s$ then $X_r\subseteq X_s$;
\item if $s$ has children $r,t$ then $(X_s\setminus (X_r\cap X_t), X_r\cap X_t)$ is a balanced skew partition of $G[X_s]$;
\item if $s\in V(T)$ has no children then either $G[X_s]$ has no balanced skew partition, or its clique number is less than $\omega(G)$,
or it is not anticonnected, or it has at most $2\omega(G)-1$ vertices.
\end{itemize}
By ``processing'' a subset $X\subseteq V(G)$, we mean doing the following:
\begin{itemize}
\item check that $|X|\ge 2\omega(G)$ (and if not, stop);
\item check that $\omega(G[X])=\omega(G)$ (and if not, stop);
\item check that $G$ is anticonnected (and if not, stop);
\item apply \ref{balancedalg} to $G[X]$ (and if there is no balanced skew partition, stop);
\item let $(A,B)$ be the output of \ref{balancedalg}; partition $A$ into two nonempty subsets $A_1,A_2$ 
such that $A_1$ is anticomplete to $A_2$, and output $A_1\cup B, A_2\cup B$.
\end{itemize}
Thus, processing a set $X$ takes time $O(n^{\max(k,6)})$, where $n=|V(G)|$ and $k = \omega(G)$.
We can construct an SP-decomposition tree by initially setting $V(T) = \{t\}$ and $X_t=V(G)$, and recursively processing each new $X_s$;
if processing $X_s$ gives us two new sets $X_r,X_t$ say we add two children $r,t$ of $s$ to $T$.
The time to find an SP-decomposition tree thus depends on the number of vertices in the tree, and that is not so easy to estimate, because
the various sets $X_t\;(t\in V(T))$ can intersect. But pellets give us a way to find a bound. We need:

\begin{thm}\label{treesize}
Let $n=|V(G)|$ and $k=\omega(G)$, and let $T, (X_t:t\in V(T))$ form an SP-decomposition tree for $G$. 
Then at most $n^{2k}-1$ vertices of $T$ have 
children, and at most $n^{2k}$ vertices have no children.
\end{thm}
\Proof
For each $t\in V(T)$, let $f(t)$ be the number of $k$-pellets included in $X_t$. We need:
\\
\\
(1) {\em $f(r) + f(t)+1\le f(s)$ for each $s\in V(T)$ with children $r,t$.}
\\
\\
Let $B = X_r\cap X_t$ and $A = X_s\setminus B$; so $(A,B)$ is a balanced skew partition of $G[X_s]$.
Every $k$-pellet of $G[X_r]$ is also a $k$-pellet of $G[X_s]$ and the same for $G[X_t]$. We claim that no $k$-pellet is counted twice.
For let $P$ be a $k$-pellet of $G[X_s]$. If it is a $k$-pellet of both $G[X_r]$ and $G[X_t]$, then $P\subseteq X_r\cap X_t = B$,
and since $P$ is anticonnected, it is a subset of some anticomponent of $G[B]$. Since $B$ has at least two anticomponents,
it follows that $\omega(G[B])>\omega(G[P])$, which is impossible since $\omega(G[P]) = k=\omega(G)$. This proves that
no $k$-pellet contributes to both $f(r), f(t)$, and so $f(r) + f(t)\le f(s)$. We need to prove
strict inequality, however; so we need to show that some $k$-pellet of $G[X_s]$ is a subset of neither of $X_r, X_t$. 
To show this, by \ref{growstable} there is a set $Y$ of cardinality $2k-1$ such that $G[Y]$ is anticonnected
and $Y$ includes a $k$-clique. By the same argument as before, $Y\not\subseteq B$, and so we may assume that $Y\cap A_1\ne \emptyset$.
Choose $a_2\in A_2$; then since $a_2$ is nonadjacent to the vertices of $Y$ in $A_1$, $Y\cup \{a_2\}$ is anticonnected and is
thus a $k$-pellet of $G[X_s]$,
and not a $k$-pellet of either of $G[X_r],G[X_t]$. This proves (1).

\bigskip
Thus, if we sum $f(s) - f(r)-f(t)-1$ over all vertices $s$ that have children $r,t$ say, then the answer is nonnegative. 
Let $q$ be the root, and $L$ the set of vertices with no children, and $I$ the set that have children. Then the sum 
can be rewritten as
$f(q) - |I| - \sum_{r\in L} f(r)$, and hence is at most $f(q)-|I|$; and so $|I|\le f(q)$. Since $f(q)< n^{2k}$, and $|L|=|I|+1$, this 
proves \ref{treesize}.~\bbox

Let us put these pieces together to prove \ref{mainalg}. We prove the following.
\begin{thm}\label{mainthm}
For all $k\ge 1$, there is an algorithm $\mathcal{P}_k$ that, with input an $n$-vertex Berge graph $G$ with $\omega(G)=k$,
outputs in time $O(n^{(k+1)^2})$ a $k$-colouring of $G$.
\end{thm}
\Proof
We proceed by induction on $k$. The result is true for $k = 1,2$, so we assume that $k\ge 3$ and the result holds
for $k-1$.  Here is the algorithm:
\begin{itemize}
\item Construct an SP-decomposition tree $T, (X_t:t\in V(T))$ for $G$. Since $T$ has at most $n^{2k}$ internal vertices by \ref{treesize},
and $2k+\max(k,6)\le (k+1)^2$, this takes time $O(n^{(k+1)^2})$.
\item For each vertex $t\in V(T)$ such that $t$ has no children, compute an optimal colouring for $G[X_t]$ as follows. Since $t$ has no children,
either $G[X_s]$ has no balanced skew partition, or its clique number is less than $k$,
or it is not anticonnected, or it has at most $2k-1$ vertices. In the first case we apply \ref{noskewalg}, in the second case
we apply $\mathcal{P}_{k-1}$, in the third we apply $\mathcal{P}_{k-1}$ to its anticomponents, and in the fourth we can find an
optimal colouring in constant time (depending only on $k$). A call to $\mathcal{P}_{k-1}$ takes time $O(n^{k^2})$ 
from the inductive hypothesis, and since $k^2\ge 7$, in each case we can find an optimal colouring in time 
$O(n^{k^2})$.
There are at most $n^{2k}$ such vertices $t$, and since $2k+k^2\le (k+1)^2$, this step takes time $O(n^{(k+1)^2})$ altogether.
\item Combine these colourings to obtain an optimal colouring of $G$ (starting from the leaves and working inwards), 
using the algorithm of section 9 (with $\mathcal{P}_{k-1}$)  at most $n^{2k}$ times. 
Each call of the algorithm from section 9 takes time $O(n^{k^2})$,
and since $2k+k^2\le (k+1)^2$, 
altogether this step takes time $O(n^{(k+1)^2})$. 
\end{itemize}
The total running time is thus $O(n^{(k+1)^2})$. 
This proves \ref{mainthm}.~\bbox

\section{Acknowledgements}
The authors would like to thanks Fr\'ed\'eric Maffray,  Nicolas Trotignon and Kristina Vu\v{s}kovi\'{c} for many helpful 
discussions.

\end{document}